\documentclass{article}
\usepackage{amssymb}
\newtheorem{definition}{\bf Definition}[section]

\newtheorem{lemma}{\bf Lemma}
\newtheorem{theorem}{\bf Theorem}

\def\Box{\square}
\newcommand{\proof}[1]{ #1 \hfill$\Box$\par}
\newcommand{\fcaption}[1]{
        \refstepcounter{figure}
           \begin{center}
            \parbox{4in}{\footnotesize\rm \baselineskip=10pt Fig.~\thefigure. #1 }
            \end{center}
        \vspace{-3pt}}

\input epsf
\epsfverbosetrue

\newcommand{\Init}[2]{{\rm stem}(#1,#2)}
\newcommand{\low}{{\rm low}}
\newcommand{\down}{{\rm Fringe}}
\newcommand{\interl}{{\rm Interlaced}}
\newcommand{\lset}{{\rm L}}

\newcommand{\out}{{\omega^+}}

\newcommand{\startlevel}{\begin{small}\begin{list}{}{\itemsep 0pt \topsep 0pt}}
\newcommand{\finishlevel}{\end{list}{}{}\end{small}}

\begin{document}

%% print out the publisher copyright heading
%\copyrightheading

%% use symbolic footnote
%\symbolfootnote

%% use normal text like skip (13pt)
%\textlineskip

%% print out titles in IJFCS format
\title{Depth-First Search and Planarity}
\author{Hubert de Fraysseix,\\
 Patrice Ossona de Mendez,\\
Pierre Rosenstiehl
\vspace{10pt}\\
{\small Centre d'Analyse et de Math\'ematiques Sociales (CNRS UMR 8557),}\\
{\small \'Ecole des Hautes \'Etudes en Sciences Sociales, 54 boulevard Raspail}\\
{\small Paris, 75006, France}}
\date{}
\maketitle
\begin{abstract}
We present a simplified version of the DFS-based Left-Right planarity testing and embedding
algorithm implemented in Pigale \cite{Taxi_pigale,Taxi_pigchap}, 
which has been considered as the fastest implemented one \cite{comp_algo03}.
We give here a full justification of the algorithm, based on a
topological properties of Tr\'emaux trees.
\end{abstract}
\section{Introduction}
It is well known since the publications \cite{HT0,HT} in 1973-74 by
J. Hopcroft and 
R.E. Tarjan that the time complexity of the problem of graph planarity
testing is linear in the number of edges. However, authors of text
books and teachers of graph theory know how hard it is to describe and
completely justify such algorithms, each one being more tricky than
the other. To display a planar embedding is usually a second hard step
\cite{Williamson-Embed}. 

A satisfactory graph planarity algorithm should be efficient, clearly
justified, and easy to understand. Only a strong mathematical insight
of the subject could probably meet these three criteria. 

R. Tarjan initiated the use of the Depth-First Search procedure (DFS) to
attack planarity testing, both DFS and planarity testing on
biconnected graphs being handled recursively. Several authors
explained how to extract an embedding during the recursion of
such an algorithm. The authors produced in the eighties the so-called
{\sc Left-Right} algorithm, a non-recursive version
avoiding the $2$-connectivity assumption. The
{\sc Left-Right} algorithm appeared to be
extremely efficient for testing planarity and embedding planar graphs
and has been recognized  as the fastest among the
implemented ones by the comparative tests
performed by graph drawing specialists \cite{comp_algo03}. 
But even if a tentative mathematical explanation
appeared in our papers on Tr\'emaux characterization of planarity
\cite{Taxi_Tremaux,Taxi_plantrem,Taxi_plandfs}, 
 it was not enough to fully meet the second and third criteria. The
code appears in the GPL-licenced software PIGALE
\cite{Taxi_pigale,Taxi_pigchap}.
 A new, simplified,
and faster version has been implemented at the occasion of this
paper, which is the one we will discuss here.

In Section~\ref{sec:Tremaux}, Tr\'emaux trees are studied as a
mathematical object {\em per se}. A rooted tree that spans a graph is a Tr\'emaux tree if 
each cotree edge is incident to two comparable vertices (with respect
to the tree order). 
It is why cotree edges are then called
back-edges. For a given Tr\'emaux tree, the following
structural concepts are defined: the low of a vertex or of an edge,
the low set of a vertex or of an edge and the fringe of an edge.
%, the partition of the
%tree edges into three classes : the block edges, the thin edges and
%the thick edges.
%
% Depth-first search is the classical way to generate a
% Tr\'emaux tree, and several DFS's may generate the same Tr\'emaux
% tree. We make the connection between these tools in
% Section~\ref{sec:Tprec} by introducing
% $T$-precedence orders. We also show that such orderings of the outgoing
% edges at each vertex of the graph and a bipartition of the back-edges
% into left- and right-attached lower incidences define particular drawings where
% crossings are easily characterized. 
% Special $T$-precedence orders are
% introduced, the $TT$-precedence orders, which are defined using the
% $\low$ function and the partition of the edges into thin, thick and
% block edges.

Section~\ref{sec:planar} is devoted to the study of planarity.
% using
%the tools introduced in the previous sections.
Planarity has been related to Tr\'emaux trees in \cite{Taxi_Tremaux} in a
characterization based on the existence of a bicoloration constrained by 
three special simple configurations. These configurations 
are here unified into a single one. 
% defined by the relative positions of $3$ or $4$ back-edges with
% respect to a Tr\'emaux tree $T$.
%  In one of the configurations
% two of the back-edges are declared to be $T$-alike, while in
% the two others, two of the back-edges are declared
% $T$-opposite. The original theorem of \cite{Taxi_Tremaux} (quoted here as
% Theorem~\ref{thm:planar})
% states that a graph with a given Tr\'emaux tree $T$
% is planar if and only if all the $T$-alike and $T$-opposite
% declarations of all its Tr\'emaux configurations are compatible, that is
% consistent with a bipartition of the set of back-edges.
% It is then proved that the bipartition constraints defined by the
% $T$-alike and $T$-opposite relations may also be expressed as an alternaltive
% bicoloration condition by introducing an arbitrary $TT$-precedence order. 
In this setting, it is proved that some additional constraints may be
imposed, which do not change the existence of a bipartition but lead a
simple planar embedding.

The planarity testing and embedding algorithm is then described in Section~\ref{sec:algo}.
It is shown that linear time-complexity is reached by implicitly
building a spanning arborescence in the graph of constraints. It shall be
noticed that the data structures used by the algorithm are almost
trivial ones, which may explain its computational efficiency and the
ease to follow the algorithm step by step on an example.

In the sequel, by a ``graph'' we mean a ``connected, loopless
multigraph'', unless we state it otherwise. 

\section{Tr\'emaux Trees}
\label{sec:Tremaux}

{\em Depth-first search} (DFS) is a fundamental graph searching
technique known since the $19$th century (see for instance Luca's
report on Tr\'emaux's work \cite{lucas82}) and popularized by Hopcroft and
Tarjan \cite{HT0,Tarjan} in the seventies. 
The structure of DFS enables efficient algorithms
for many other graph problems \cite{even_book}. Performing a DFS on a 
graph defines a spanning tree with specific properties (also known as
a {\em Tr\'emaux tree}) and an embedding of it as a rooted planar
tree, the edges going out of a vertex being circularly ordered
according to the discover order of the DFS).

%\subsection{Tr\'emaux tree partial order}

A rooted spanning tree $T$ of a graph $G$
defines a partition of the edge set of $G$ into two classes, the set
of {\em tree edges} $E(T)$ and the set of {\em cotree-edges} $E(G)\setminus E(T)$.
It also defines a partial order $\preceq$ on $V(G)$: $x\preceq y$ if the
tree path linking $y$ to the root of $T$ includes $x$. The rooted tree $T$ is
a {\em Tr\'emaux tree} if every cotree edge is incident to two comparable vertices
(with respect to $\preceq$). A Tr\'emaux tree $T$ defines an orientation of the edges
of the graph: an edge $\{x,y\}$ (with $x\prec y$) is oriented from $x$
to $y$ (upwards) if it is
a tree edge and from $y$ to $x$ (downwards) if it is a cotree edge. Cotree edges
of a Tr\'emaux tree are called {\em back edges}.
We will denote by $\out(v)$ the set of the edges incident to a vertex $v$ and going out of $v$.
When $T$ is a Tr\'emaux tree, the partial order
$\preceq$ is extended to $V(G)\cup E(G)$ (or to $G$ for short) as follows: 
for any edge $e=(x,y)$ oriented
from $x$ to $y$, put $x\prec e$ and if $x\prec y$ 
(that is: if $e$ is a tree-edge) also put $e\prec y$.
\begin{figure}[ht]
\begin{center}
 \epsfxsize=0.5\textwidth
\epsffile{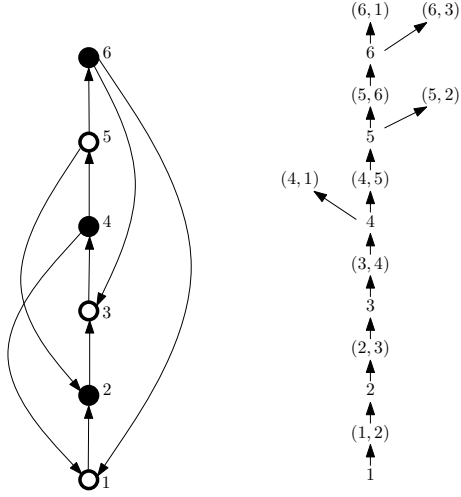}
\end{center}
\caption{The partial order $\prec$ defined by a Tr\'emaux tree of
  $K_{3,3}$ (cover relation correspond to bottom-up arrows as usual).}
\label{fig:prec}
\end{figure}

Notice that in the partial order $\prec$, all elements $\alpha$ and $\beta$ of $G$ have 
an unique greatest lower bound (meet) $\alpha\wedge\beta$.
 Moreover, the maximal chains of $\prec$ all have the
same structure: they begin with the root vertex of $T$, alternate
between vertices and edges, and include at most one back-edge (which
is then the maximum of the chain).

When $\alpha\preceq\beta$, the unique chain (of $\prec$) with minimum
$\alpha$ and maximum $\beta$ which is maximal (with respect to set-inclusion)
is denoted by $[\alpha\,;\beta]$.

It will be helpful to
introduce a notation for the minimal element of the interval $]\alpha\,;\beta]$:
\begin{definition}
\label{def:init}
For $x\prec e$, where $x\in V(G)$ and $e\in E(G)$ we define
\begin{equation}
\Init{x}{e}=\min\ ]x\,;e]
\end{equation}
This means that $f$ is the first edge in the unique chain of $\prec$ with minimum
$x$ and maximum $e$.
\end{definition}
Notice that in this definition, as in the remaining of the paper,
intervals and ``$\min$'' will always be related to the partial order
$\prec$ defined by the considered Tr\'emaux tree.

%\subsection{Related Definitions}

\begin{definition}
\label{def:low}
The function $\low:G\rightarrow V(G)$ is defined by
\[
\low(\alpha)=\min\ \bigl(\{\alpha\}\cup\{v\in V(G): \exists  
(u,v)\in E(G)\setminus E(T),\ (u,v)\succeq\alpha\,\}\bigr).
\]
Notice that this function is well defined and that $\low(\alpha)\preceq
\alpha$. 
\end{definition}

In the literature, this function is usually only defined on
$V(G)$. Our extension is such that for any edge $e=(x,y)$:
$$\low(e)=\left\{\begin{array}{ll}
\min(x,\low(y)) &\mbox{if $e$ is a tree edge,}\\
y &\mbox{if $e$ is a back-edge.}
\end{array}\right.$$

%We shall need more information on how a vertex is linked to its $\low$.

\begin{definition}
The {\em fringe} $\down(e)$ of an edge $e=(x,y)$ is defined by:
\[
\down(e)=\{f\in E(G)\setminus E(T): f\succeq e\mbox{ and
}\low(f)\prec x\,\}
\]
\end{definition}

\begin{definition}
The {\em low set} $\lset(\alpha)$ of $\alpha\in V(G)\cup E(G)$ is defined by:
\[
\lset(\alpha)=\{f\in E(G)\setminus E(T): f\succeq\alpha\mbox{ and }\low(f)=\low(\alpha)\,\}
\]
Notice that $\lset(e)=\{e\}$ if $e$ is a back-edge.
\end{definition}

{\flushleft
\begin{tabular}{ll}
\parbox[b]{0.6\textwidth}{
\begin{definition}
 Given a tree edge $e\in E(T)$, 
we will call $e=(x,y)$:

\begin{itemize}
\item 
a {\em block edge} 

if $x\preceq\low(y)$, that is if $low(e)=x$
(this means $e$ is either an isthmus or the minimum edge of a block of
the graph);
%notice that there is exactly one such edge per block of the graph);
\item a {\em thin edge} 

if $\low(y)\prec x$ and there exists no
  back-edge $(u,v)$ with $u\succeq y$ and $\low(y)\prec u\prec x$;
\item a {\em thick edge} 

if $\low(y)\prec x$ and there exists a
  back-edge $(u,v)$ with $u\succeq y$ and $\low(y)\prec u\prec x$.
\end{itemize}
\end{definition}
}&
\parbox[b]{0.4\textwidth}{
 \epsfxsize=0.35\textwidth
 \epsffile{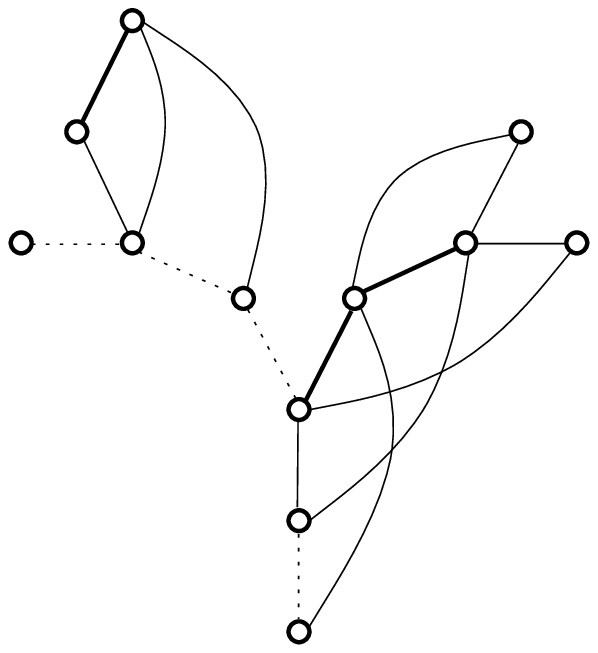}
\fcaption{
\parbox[t]{0.35\textwidth}{Block, thin and thick edges.
Block edges are dotted, thin edges
are light and thick edges are fat.
\label{fig:btt}}}
\bigskip

}
\end{tabular}
}
%\end{definition}
%\bigskip

\begin{definition}
A {\em $TT$-precedence order} $\prec^\star$ is a partial 
order on $E(G)$ such that, for any $v\in V(G)$ and any $e,f\in\out(v)$:
\begin{itemize}
\item if $\low(e)\prec\low(f)$ then $e\prec^\star f$,
\item if $\low(e)=\low(f)$, $f$ is a thick tree edge but $e$ is not,
  then $e\prec^\star f$.
\end{itemize}
\end{definition}
\clearpage
\section{Tr\'emaux Trees and Planarity}
\label{sec:planar}
%\subsection{Tr\'emaux Trees and Planarity}
%\label{sec:TrEmauxTreesAndPlanarity}

Planarity has been related to Tr\'emaux trees by de Fraysseix and
Rosenstiehl in a series of articles
\cite{Taxi_Tremaux,Taxi_plantrem,Taxi_plandfs}. 
One of these characterizations is based on a the existence of a
special bipartition of the low
angles of the back-edges into left ones and right ones. The
constraints that the bipartition has to fulfill is encoded into two
relations, namely the $T$-alike and $T$-opposite relations.

\begin{figure}[ht]
   \begin{center}
     \fbox{\epsfxsize=0.5\textwidth\epsffile{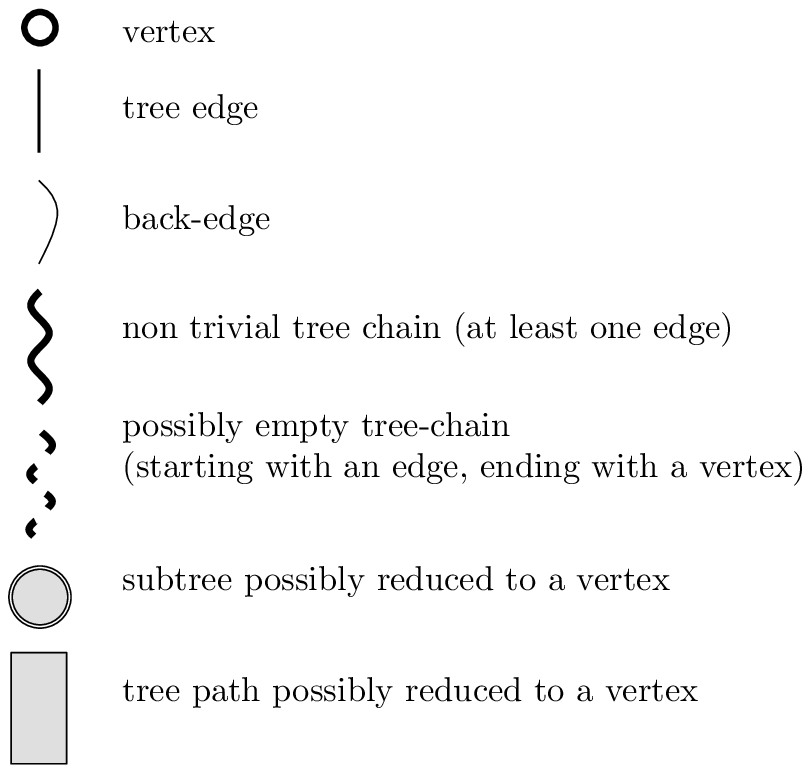}}
   \end{center}
   \caption{Symbols used in the figures}
   \label{fig:legend}
\end{figure}

\begin{figure}[ht]
   \begin{center}
     \begin{tabular}{ccc}
       \epsfxsize=71.2pt
       \epsffile{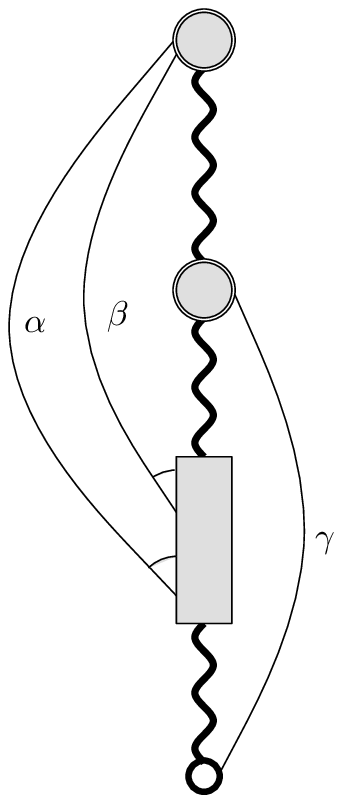}&
       \epsfxsize=87pt
       \epsffile{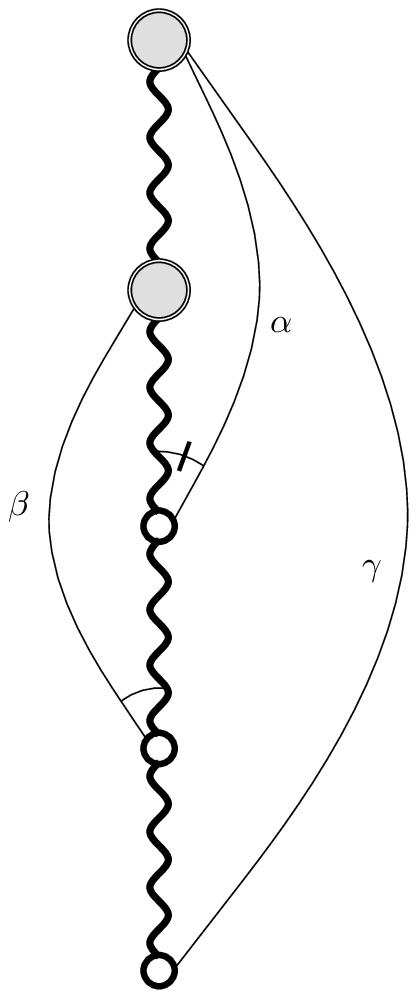}&
       \epsfxsize=101.6pt
       \epsffile{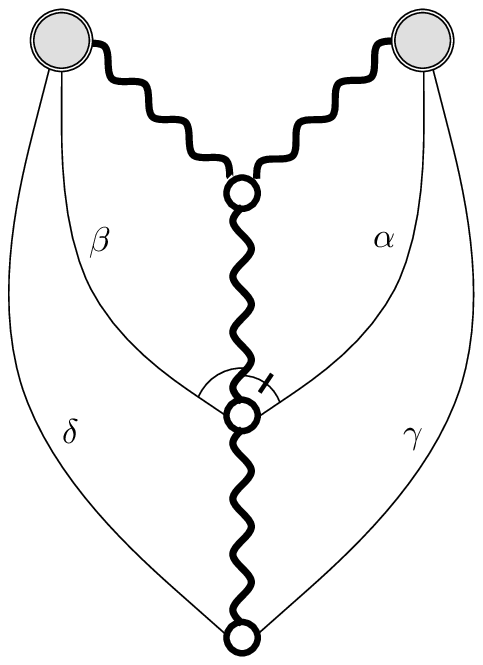}\\
     case (i)&case (ii)&case (iii)\\
     {\small\rm $\alpha$ and $\beta$ are $T$-alike}&
     {\small\rm $\alpha$ and $\beta$ are $T$-opposite}&
     {\small\rm $\alpha$ and $\beta$ are $T$-opposite}
   \end{tabular}
\end{center}
\caption{Definition of $T$-alike and $T$-opposite relations\label{fig:tao}}
\end{figure}

We don't give here the formal definition of $T$-alike and $T$-opposite relations in terms 
of $\prec$, but simply recall the characterization given in \cite{Taxi_Tremaux}:
\begin{theorem}
\label{thm:planar}
Let $G$ be a graph with Tr\'emaux tree $T$.
Then $G$ is planar if and only if there exists a partition
of the back-edges of $G$ into two classes so that any two edges
belong to a same class if they are $T$-alike and any two edges belong to
different classes if they are $T$-opposite. 
\end{theorem}

Instead of working with this characterization, we introduce an
equivalent characterization based on a single configuration.
\begin{definition}
Let $v$ be a vertex and let $e_1,e_2\in\out(v)$.

The {\em interlace set} $\interl(e_1,e_2)$ is defined by:
\[
\interl(e_1,e_2)=\{f\in \down(e_1): f\succ\low(e_2)\,\}
\]
\end{definition}

\begin{definition}
Given a graph $G$ and a Tr\'emaux tree $T$ of $G$,
a coloring $\lambda:E(G)\setminus E(T)\rightarrow\{-1,1\}$ is an 
{\em F-coloring} if, for every vertex $v$ and any edges $e_1,e_2\in\out(v)$,
$\interl(e_1,e_2)$ and $\interl(e_2,e_1)$ are
monochromatic and colored differently.
\end{definition}

It is easily checked that a coloring is an $F$-coloring if and only if
two $T$-alike back-edges are colored the same and two $T$-opposite
back-edges are colored differently.

In a planar drawing, an $F$-coloring is defined by the partition of
the back-edges on two sets, the edges  $f$ having their low incidence on the
left (resp. the right) of the tree edge $\Init{\low(f)}{f}$.
The following lemma is straightforward and does not deserve a proof:

\begin{lemma}
\label{lem:planar}
Let $G$ be a planar graph with Tr\'emaux tree $T$.
Then $G$ has an $F$-coloring.
\hfill$\Box$
\end{lemma}

\begin{definition}
An F-coloring $\lambda:E(G)\setminus E(T)\rightarrow\{-1,1\}$ is {\em
  strong} if, for any $v\in V(G)$ the low set $\lset(v)$ is
monochromatic.
\end{definition}

\begin{lemma}
If $G$ has an F-coloring then it has a strong F-coloring.
\end{lemma}
\proof{The addition of the constraints that the sets $\lset(v)$ are
  monochromatic may not lead to a contradiction, as the only
  constraints involving $f_1\in\lset(v)$ would also involve any
  $f_2\in\lset(v)$ and would require that $f_1$ and $f_2$ actually
  have the same color.}

%\subsection{Computing an Embedding from a Strong $F$-coloring}

>From the ``low angles bicoloration'' $\lambda$ we define a ``high angles bicoloration'' 
$\widehat\lambda$ on the whole edge set of the graph (on both tree edges and back-edges).

\begin{definition}
Let $\lambda:E(G)\setminus E(T)\rightarrow\{-1,1\}$ be a strong
F-coloring. We define the coloring
$\widehat\lambda:E(G)\rightarrow\{-1,1\}$ by:
\[\widehat\lambda(e)=
\left\{\begin{array}{ll}
  \lambda(e),&\mbox{ if $e$ is a back-edge}\\
  \lambda(f),&\mbox{ if $e$ is a tree edge and }f\in\down(e)\mbox{ with
    maximal }\low(f)\\
\end{array}\right.
\]
\end{definition}

\begin{lemma}
\label{lem:embed}
 Let $G$ be a graph, let $T$ be a Tr\'emaux tree
 of $G$, let $\lambda$ be a strong
$F$-coloring and let $\widehat\lambda$ be the associated mapping.
We define the circular order of the edges at a vertex $v$ as follows:

Let $e_1\succ^\star e_2\succ^\star\dots\succ^\star e_p$ be the
edges in $\out(v)$ with $\widehat\lambda(e_i)=-1$ and let
$e_{p+1}\prec^\star e_{p+2}\prec^\star\dots\prec^\star e_q$ be the
edges in $\out(v)$ with $\widehat\lambda(e_i)=1$.

In the circular
order around $v$ one finds the incoming tree edge (if $v\neq r$) and then
$L_1,e_1,R_1,L_2,e_2,R_2,\dots,L_q,e_q,R_q$ where $L_i$ (resp. $R_i$)
is the set of incoming cotree edges $f=(x,v)$ such that
$\lambda(f)=-1$ (resp. $\lambda(f)=1$) and the tree-path linking
$r$ to $x$ includes $e_i$.
For $e_i,e_j\in L_k$ (resp. for $e_i,e_j\in R_k$), 
one finds $e_i$ before (resp. after) $e_j$ in the circular order if
$\Init{e_i\wedge,e_j}{e_i}\prec^\star\Init{e_i\wedge e_j}{e_j}$.

Then these circular orders define a planar embedding of $G$.
\end{lemma}

\proof{In a drawing where the tree edges cross no other edges, 
only two kind of crossings could occur:
\begin{center}
\begin{tabular}{cc}
\parbox[b]{0.5\textwidth}{\epsffile{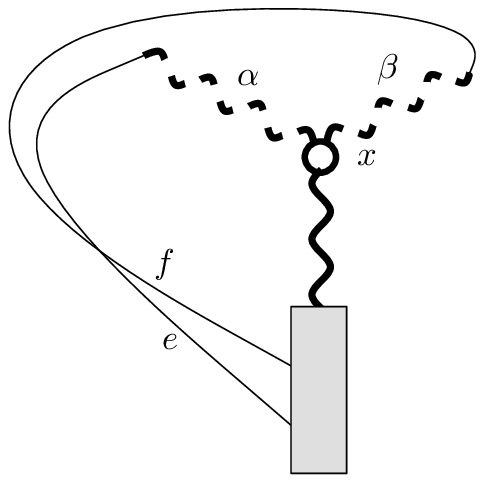}}&
\parbox[b]{0.5\textwidth}{\epsffile{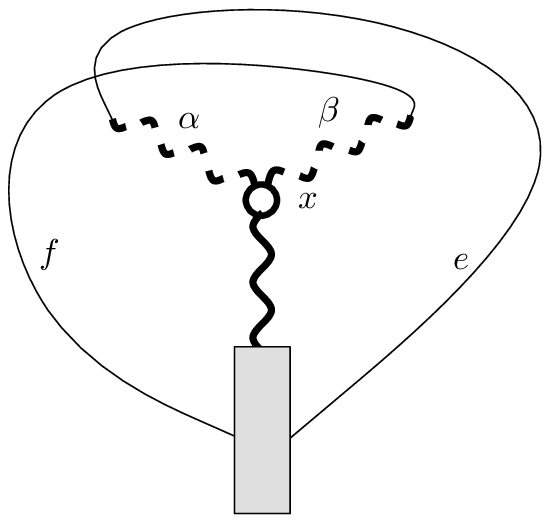}}
\end{tabular}
\end{center}
%We shall assume that $G$ contains exactly one block edge
%incident to the root (by considering each block independently).
Let $x=e\wedge f, y=\low(f), z=\low(e), \alpha=\Init{x}{e}$ and
$\beta=\Init{x}{f}$.
Without loss of generality we assume  $\lambda(f)=-1$.
\begin{itemize}
\item 
The first case corresponds to two back-edges $e,f$ with
$\lambda(e)=\lambda(f)$ and $\low(e)\leq\low(f)$.

If $\low(\alpha)=\low(\beta)=\low(e)=\low(f)=u$, then the
contradiction arises from the definition of the circular order at $u$.
Otherwise, as $f\in\interl(\beta,\alpha)$, $f$ is colored the same way as the
back-edge in $\down(\beta)$ with maximal $\low$-value. Hence
$\widehat\lambda(\beta)=\lambda(f)=-1$. According to the definition of
the circular order at $x$, $\widehat\lambda(\alpha))=-1$ and
$\alpha\succ^\star\beta$. According to the thin-thick precedence of
$\prec^\star$ it follows that the edge $e'\in\down(\alpha)$ with
maximal $\low$-value is such that $\low(e')\succ\low(\beta)$, that is:
$e'\in\interl(\alpha,\beta)$. 
Thus $\widehat\lambda(\alpha)=\lambda(e')\neq\lambda(f)$, a contradiction.

\item 
The second case corresponds to two back-edges $e,f$ with
$\lambda(e)\neq\lambda(f)$. By symmetry we may assume
$\low(\alpha)\leq\low(\beta)$ and $\lambda(e)=1$.

Then $\low(\beta)\succ\low(x)$ as $\lambda$ is a strong $F$-coloring.
Thus $f$ is colored the same way as the
back-edge in $\down(\beta)$ with maximal $\low$-value and 
$\widehat\lambda(\beta)=\lambda(f)=-1$. According to the definition of
the circular order at $x$, $\widehat\lambda(\alpha))=-1$ and
$\alpha\succ^\star\beta$. It follows that $e\in\lset(x)$, a contradiction.
\end{itemize}
}

 \begin{theorem}
 \label{thm:planar2}
 Let $G$ be a graph, let $T$ be a Tr\'emaux tree
 of $G$. The following conditions are equivalent:
\begin{enumerate}
\item[(i)] $G$ is planar,
% \item[(ii)] there exists a bipartition of $E(G)\setminus E(T)$ such
%   that any two $T$-alike edges belong to a same part and any two
%   $T$-opposite edges belong to different parts,
\item[(ii)] $G$ admits an F-coloring,
\item[(iii)] $G$ admits a strong $F$-coloring.
\end{enumerate}

Moreover, if $G$ is planar, any strong F-coloring $\lambda$ defines a
planar embedding of $G$ in which a back-edge $e$ has its lowest incidence
to the left of the tree if $\lambda(e)=-1$ and to the right of the
tree if $\lambda(e)=1$·
 \end{theorem}
 \proof{This is a direct consequence of
   the previous lemmas.}

% \subsection{Embedding}
% \begin{definition}
% A {\em good coloring} is a red-blue coloring of the back-edges, such
% that any two $T$-alike edges get the same color and any $T$-opposite
% edges get different colors.

% A good coloring is {\em strong} if, for every vertex $v$, all the
% edges $e=(x,\low(v))$ with $x\succeq v$ get the same color.
% \end{definition}

% \begin{lemma}
% Let $v$ be a vertex and let $e,f\succeq v$ be two back-edges with
%  lower incidence $low(v)$. Then either $e$ and $f$ may have
%  (independently) any color in a good coloration, or $e$ and $f$ always
%  have the same color in a good coloration.

% Consequently, if $G$ has a good coloration, it has a strong good coloration.
% \end{lemma}

\section{The Planarity Testing Algorithm}
\label{sec:algo}
\subsection{Outline}
% The outline of the planarity testing algorithm is the following:

% \begin{itemize}
% \item {\bf do} a preliminary DFS on $G$, compute the $\low$ and the status of the edges
%   (block/thin/thick).
% \item {\bf do} a sort of the edges according to a $TT$-precedence order.
% \item {\bf do} a DFS on $G$ steered by the $TT$-precedence order
%   (using the sort performed at the previous step) and maintain the
%   bicoloration constraints.
% {\bf If} the bicoloration constraints did not lead to a
%   contradiction
% {\bf then} output that $G$ is planar
% {\bf else} output that $G$ is not planar.
% \end{itemize}

Let $G$ be a graph of size $m$. The three steps
are performed in $O(m)$-time.
The first step is composed of a preliminary DFS on $G$ and the
computation of the $\low$ function and the status of the edges (block/thin/thick).
The second step is the computation of a $TT$-precedence
order, which may be efficiently performed using a bucket
sort. We now examine the last step of the algorithm, which tests the
planarity of the graph.

%  and only the third step requires a bit of
% attention.

% The DFS considered in the third step 
% discovers between two backtracks (on tree on cotree edges) some walks
% $W_e$, which form a partition of the edge set of $G$.
% Theses walks are either reduced to a cotree edge or  
% start with some tree edge $e$, follow (as long as possible) reference
% tree edges and end with a reference back-edge (if the tree terminal has degree at least $2$).

% Bicoloration constraints related to these walks are stored into {\em Double Stack
%   Systems} (DSS), each DSS being a stack of {\em Double Stacks} (DS)
% (i.e. pairs of stacks). A DS contains a subset of back-edges which
% are dependent (either constrained to have a same color if they belong
% to the same stack of the pair,  or constrained to have different
% colors if they belongs to different stacks of the pair). 
% The DS are stacked upon one another, meaning that no constraints exist
% between the back-edges in one DS and a back-edge in another at the
% considered moment of the algorithm.

% % A DSS is created for each walk, including a single DS
% % with the back-edge of the walk (if any).

% \subsection{New}

We shall consider some data structure ${\rm CS}$ responsible for maintaining a set
of bicoloration constraints on a set of back-edges. We assign
to each edge of the graph $e$ such a data structure ${\rm CS}(e)$. These structures are
initialized as follows: ${\rm CS}(e)$ is empty if $e$ is a tree edge
and includes $e$ (with no bicoloration constraints) if $e$ is a
back-edge. We say that all the back-edges have been {\em processed}
and that the tree edges are still {\em unprocessed}.

\begin{itemize}
\item While there exists a vertex $v$, different from the root, such that 
all the edges in $\out(v)$ have been processed.
Let $e=(u,v)$ be the tree edge entering $v$. 
Let $e_1\prec^\star e_2\prec^\star\dots\prec^\star e_k$
be the edges in $\out(v)$ ($k\geq 1$). We do the following:
\begin{itemize}
\item Initialize ${\rm CS}(e)$ with ${\rm CS}(e_1)$.
\item For $i:2\rightarrow k$, {\em Merge} ${\rm CS}(e_i)$
  into ${\rm CS}(e)$, that is: add to ${\rm CS}(e)$ the edges in
  ${\rm CS}(e_i)$ and add the $F$-coloring constraints 
  corresponding to the pairs of edges $e_j,e_i$ with
  $j<i$ (notice that all the concerned back-edges belong to ${\rm CS}(e)$).
  If some constraint may not be satisfied, the
  graph is declared non-planar.
\item Remove from all of the ${\rm CS}(e)$ every back-edge with
  lower incidence $u$.
\item We declare that edge $e$ has been processed.
\end{itemize}

\item As all the edges have been processed, we declare that the graph is planar.
\end{itemize}

\subsection{Data-structure, Complexity and Embedding Computation}

First
notice that the processing order of the tree edges is simple to
compute, using either a ``topological sort technique'' (by maintaining
the unprocessed outdegree of vertices) or by following the backtrack
order of a DFS steered by the TT-order. Also, we should notice that
the structure ${\rm CS}(f)$ is only used (after its computation) 
when computing ${\rm CS}(e)$ where $f$ is the predecessor of $f$ in
$\prec$. It follows that we may ``destroy'' ${\rm CS}(f)$ when
computing  ${\rm CS}(e)$ without any risks. 

%We describe here the data structure ${\rm CS}$ which allows to obtain
%to linear time implementation of the algorithm.

The data structure for ${\rm CS}(e)$ is a stack of double top-to-bottom linked stacks. This
means that each element of ${\rm CS}(e)$ is a pair $({\rm S}^0,{\rm S^1})$
of stacks. Each such pair correspond to a complete
bipartite constraint graph: all the edges in ${\rm S}^i$
($i\in\{0,1\}$) have to be
colored the same and they have to be colored differently from the
back-edges in ${\rm S}^j$ (for $j\neq i$). These will be all the
constraints encoded by the ${\rm CS}(e)$ (hence no constraint exists
between back-edges belonging to different stack pairs).

Some additional constraints will allow to get amortized constant time
operations: Let $(S_1^0,S_1^1),(S_2^0,S_2^1),\dots,(S_k^0,S_k^1)$ be
the pairs stacked in some ${\rm CS}(e)$. Then:
\begin{itemize}
\item The top back-edges of $S_i^0$ of $S_i^1$ have both a lower
  incidence which is strictly smaller than the lower incidences of the
  bottoms of $S_{i+1}^0$ and $S_{i+1}^1$.
\item the lower incidences of the back-edges belonging to some stack
  $S_i^\alpha$ are in non-decreasing order.
\end{itemize}

Notice that, for any tree-edge $e=(u,v)$, ${\rm CS}(e)$ will exactly
include the back-edges $f=(x,y)$ such that $x\succ v$ and $u\succ
y$. Hence the set of edges in ${\rm CS}(e)$ is $\down(e)$.

Now consider the operations performed by the algorithm on the
structures ${\rm CS}(e)$:
\begin{itemize}
\item Merge of ${\rm CS}(e_i)$ into ${\rm CS}(e)$: 
the structure ${\rm CS}(e)$ then contains the back-edges in
$X(e_i)=\bigcup_{j<i}\down(e_j)$ and their bicoloration
constraints. By induction, the graph of these constraints is a
disjoint union of complete bipartite graphs encoded by the pairs of
stacks. The $F$-coloring conditions now express as:
\begin{itemize}
\item all the back-edges in ${\rm CS}(e)$ (i.e. in $X(e_i)$)
 which have a lower incidence strictly greater than $\low(e_i)$ should have the same color,
\item all the back-edges in ${\rm CS}(e_i)$ which have a lower
  incidence strictly greater than $\low(v)$ should have the same
  color,
\item the two above sets of edges should be colored differently.
\end{itemize}
These constraints are added as follows:
\begin{itemize}
\item let $a$ be the smallest lower incidence of the back-edges in
  ${\rm CS}(e)$ (according to monotonicity we just have to look at the
  bottom back-edges of the stacks of the bottom pair). Then
  $a=\low(v)$.
\item let $b$ be the smallest lower incidence of the back-edges in
  ${\rm CS}(e_i)$ (as above, this is computed in constant time). 
  Then $b=\low(e_i)$.
\item if ${\rm
    CS}(e)=((S_1^0,S_1^1),(S_2^0,S_2^1),\dots,(S_k^0,S_k^1))$, let $j$
  be the biggest integer $\leq k$ such that none of the back-edges in
  the top of $S_j^0$ and
  $S_j^1$ have a lower incidence greater than $b$. Then it should be checked that for
  $j+1<j'\leq k$, one of $S_{j'}^0$ and $S_{j'}^1$ is empty (for
  otherwise, we have found a contradiction in the constraints proving that
  $G$ is not planar). By flipping pairs if necessary we assume
  $S_{j'}^1$ is empty. Moreover, only one of $S_{j+1}^0$ and
  $S_{j+1}^1$ contains back-edges whose lower incidence is greater than
  $b$ (for otherwise we have a contradiction proving that $G$ is not
  planar). Up to a flipping of the pair of stack, we may assume this
  is $S_j^0$. Then we fuse $S_j^0,\dots,S_k^0$ into $S_j^0$ and
 ${\rm
    CS}(e)=((S_1^0,S_1^1),(S_2^0,S_2^1),\dots,(S_j^0,S_j^1))$.
\item similarly, if ${\rm
    CS}(e_i)=((T_1^0,T_1^1),(T_2^0,T_2^1),\dots,(T_p^0,T_p^1))$, let $q=2$
  if $T_1^0$ or $T_1^1$ has a bottom edge whose lower incidence is $a$ and
  let $q=1$ otherwise.
  Then it should be checked that for
  $q\leq q'\leq p$, one of $S_{q'}^0$ and $S_{q'}^1$ is empty (for
  otherwise, we have found a contradiction in the constraints proving that
  $G$ is not planar). By flipping pairs if necessary we assume
  $S_{q'}^0$ is empty. It should be noticed that if $q=2$, 
  then one of $T_1^0$ and $T_1^1$ is empty. By flipping the pair if
  necessary, we may assume this is
  $T_1^0$. Then $T_1^1$ contains
  exactly the back-edges in $\lset(e_i)$.
  %(this is implied by the fact
  %that requesting that an $F$-coloring is strong may never lead to a
  %contradiction). 

  Then we fuse $T_q^1,\dots,T_p^1$ into $S_j^1$ and, if $q=2$, we add
   to the one of $S_1^0$ and $S_1^1$ which is not
  empty the edges from $T_1^0$.
\end{itemize}
As the number of pairs of stacks decreases at each ``fuse'', as no new
pair of stacks is created and as the initial number of pairs of stacks
is $O(m)$, the global time spent in this step by the algorithm is $O(m)$. 
\item Deletion of the back-edges whose lower incidence is $u$: by
  monotonicity, we only have to check the top of the stacks in the top
  pair of stacks. As every back-edge will be deleted exactly once by
  the algorithm, the global time spent for deletions is $O(m)$.
\end{itemize}

As the stacks used by the algorithm have their
elements linked in a top-bottom manner, it is easy to keep, while
``deleting'' back-edges, the link between edges that have to get the
same $\lambda$-value and to add a special type of links between
some edges which have to get different $\lambda$-value. This way, a
spanning forest of the constraint graph is maintained with a constant
time cost per deletion, which allows to propagate $\lambda$-values
after the planarity testing phase. A
$\lambda$-value for each back-edge being computed, the embedding
follows and is easily computed in linear time.

%We also use a weak monotonicity property which will be fundamental to
% get a linear time algorithm.

% \begin{lemma}
% \label{lem:helpalg}
% Let $G$ be a graph, let $T$ be a Tr\'emaux tree
% of $G$ and let $\disc$ be the discover order of a DFS steered by a
% $TT$-precedence order $\preceq_0^\star$ of $G$.

% Let $e,f$ be back-edges with $\disc(f)<\disc(e)$.
% Assume $E(G)\setminus E(T)$ is bi-colored in such a way that for any non-reference
% edge $g=(a,b)\preceq e$:
% \begin{itemize}
% \item the set of $\down(g)$ is monochromatic,
% \item the color of any back-edge $h=(u,v)$ with $\disc(h)<\disc(g)$ and $\low(g)\prec v\prec
%   a\prec u$ is different from the color of the elements of
%   $\down(g)$.
% \end{itemize}

% If $e$ and $f$ have the same color then either
% $\low(e)\not\prec\low(f)$ or $\low(e)=\low(e\wedge f)$ and
% $\Init{e\wedge f}{e}$ is a thick edge.
% \end{lemma}
% \proof{
% Assume $\low(e)\prec\low(f)$. Let $x=e\wedge f, e'=\Init{x}{e}$ and
% $f'=\Init{x}{f}$. As $\disc(f)<\disc(e)$ we have $\disc(f')<\disc(e')$, hence
% $e'$ is not a reference edge and
% $\low(x)\preceq\low(f')\preceq\low(e')\preceq\low(e)$.
% If $e\in\down(e')$ we are led to a contradiction, as then $e$ and $f$
% have to be colored differently. Hence $\low(e)=\low(x)$. Because of
% the precedence of thin edges over thick edges, it follows that
% $e'$ is a thick edges (as $f'$ is a thick edge).
% }
%\subsection{Computing a Planar Embedding}
\clearpage
\subsection{Example}
\begin{figure}[h]
\begin{center}
\begin{tabular}{cc}
\epsffile{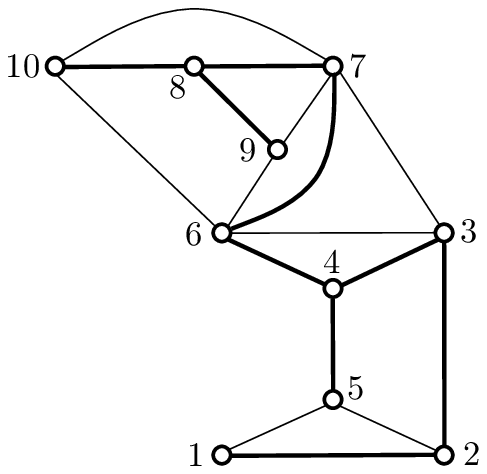}&\epsffile{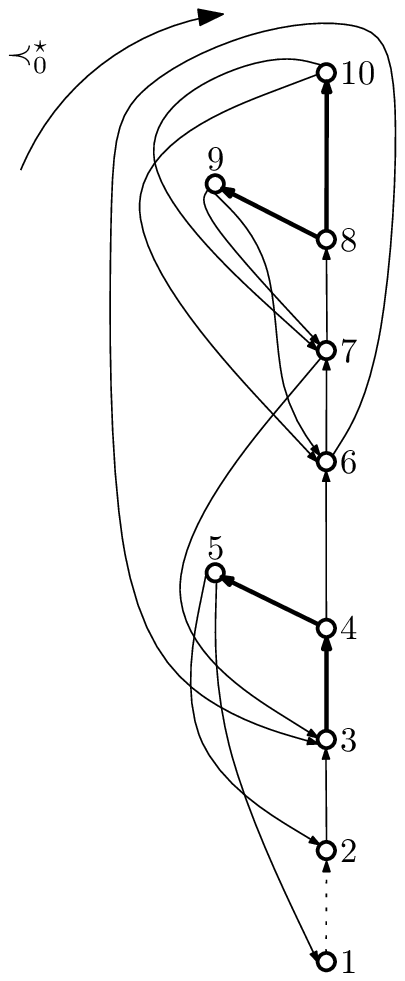}
\end{tabular}
\end{center}
\caption{A sample graph with a Tr\'emaux tree $T$ (on the left). The
  block/thin/thick partition and a $TT$-precedence order $\prec^\star$
  (on the right).
\label{fig:sample}
}
\end{figure}

% \begin{figure}[h]
% \begin{center}
% \epsffile{planar0.eps}
% \end{center}
% \caption{A sample planar graph $G$}
% \end{figure}
We consider a sample planar graph with an arbitrary Tr\'emaux tree
(see Fig.~\ref{fig:sample}).
The right of the figure display the information computed by the two
first steps of the algorithm: the status of the edges
(block/thin/thick) and a $TT$-precedence order $\prec^\star$,
represented here as a circular 
order of the outgoing edges.

The algorithm first computes the ${\rm CS}$'s of the back-edges. Then it
proceeds by iterating on the vertices whose outgoing
edges are all processed. At each such vertex $v$, the ${\rm CS}$ of
the incoming tree edge is computed, using a Merge and possibly a
Deletion step (denoted by $\rightarrow$). In the following table, the
${\rm CS}$'s are represented as sequences of double stacks.
\newcommand{\DS}[1]{\begin{array}{|c|c|}#1\\\hline\end{array}}
\clearpage
\[
\begin{array}{lclclclclcl}
v=9&:&{\rm CS}((8,9))&=&
\DS{\scriptstyle (9,6)&~~~~~}
&,&
\DS{\scriptstyle (9,7)&~~~~~}\\
\\
v=10&:&{\rm CS}((8,10))&=&
\DS{\scriptstyle (10,6)&~~~~~}
&,&
\DS{\scriptstyle (10,7)&~~~~~}\\
\\
v=5&:&{\rm CS}((4,5))&=&
\DS{\scriptstyle (5,1)&~~~~~}
&,&
\DS{\scriptstyle (5,2)&~~~~~}\\
\\
v=8&:&{\rm CS}((7,8))&=&
\DS{\scriptstyle (10,6)&~~~~~\\
\scriptstyle (9,6)&}
&,&
\DS{\scriptstyle (9,7)&\scriptstyle (10,7)}\\
\\
&&&\rightarrow&\DS{\scriptstyle (10,6)&~~~~~\\
\scriptstyle (9,6)&}&&\quad\mbox{(Deletion)}\\% of $(9,7)$ and $(10,7)$)}\\
\\
v=7&:&{\rm CS}((6,7))&=&
\DS{\scriptstyle (7,3)&~~~~~}
&,&
\DS{\scriptstyle (10,6)&~~~~~\\
\scriptstyle (9,6)&}\\
\\
&&&\rightarrow&\DS{\scriptstyle (7,3)&~~~~~}
&&\quad\mbox{(Deletion)}% of $(9,6)$ and $(10,6)$)}
\\
\\
v=6&:&{\rm CS}((4,6))&=&
\DS{\scriptstyle (6,3)&~~~~~\\
\scriptstyle (7,3)&}\\
\\
v=4&:&{\rm CS}((3,4))&=&
\DS{\scriptstyle (5,1)&~~~~~}
&,&
\DS{\scriptstyle (5,2)&~~~~~}
&,&
\DS{\scriptstyle (6,3)&~~~~~\\
\scriptstyle (7,3)&}\\
\\
&&&\rightarrow&
\DS{\scriptstyle (5,1)&~~~~~}
&,&
\DS{\scriptstyle (5,2)&~~~~~}
&&\quad\mbox{(Deletion)}% of $(6,3)$ and $(7,3)$)}
\\
\\
v=3&:&{\rm CS}((2,3))&=&
\DS{\scriptstyle (5,1)&~~~~~}
&,&
\DS{\scriptstyle (5,2)&~~~~~}\\
\\
&&&\rightarrow&
\DS{\scriptstyle (5,1)&~~~~~}
&&\quad\mbox{(Deletion)}% of $(5,2)$)}
\\
\\
v=2&:&{\rm CS}((1,2))&=&
\DS{\scriptstyle (5,1)&~~~~~}\\
\\
&&&\rightarrow&\emptyset
&&\quad\mbox{(Deletion)}% of $(5,1)$)}
\end{array}
\]

\section*{References}
%\bibliographystyle{unsrt}
%\bibliography{bib,biblio}

\end{document}